\documentclass[reqno]{amsart}
\usepackage{hyperref}
\usepackage{bm}
\usepackage{amssymb}

\begin{document}
\title[\hfil Concrete realizations]
{ Mikusi\'nski's Operational Calculus with Algebraic Foundations and Applications to Bessel Functions }

\author[Gabriel Bengochea, Gabriel L\'opez G.]
{Gabriel Bengochea, Gabriel L\'opez G.}

\address{ Gabriel Bengochea, Universidad  Aut\'onoma de la Ciudad de M\'exico. \newline
Gabriel L\'opez,  Universidad Aut\'onoma Metropolitana , M\'exico D.F, M\'exico}
\email{gabriel.bengochea@uacm.edu.mx, gabl@xanum.uam.mx}

\thanks{ } 

\begin{abstract}
We construct an operational calculus supported on the algebraic operational calculus introduced by Bengochea and Verde in \cite{BV}. With this operational calculus we study the solution of certain Bessel type equations.
\end{abstract}

\maketitle \numberwithin{equation}{section}
\newtheorem{theorem}{Theorem}[section]
\newtheorem{lemma}[theorem]{Lemma}
\newtheorem{definition}{Definition}[section]
\newtheorem{remark}{Remark}
\newtheorem{example}{Example}

\section{Introduction}
Many interesting ordinary differential equations can be solved using Mikusiński's\break operational calculus and recently, several approaches to the construction of different Calculi have been developed  (see for instance \cite{B}, \cite{D}, \cite{Mi}, \cite{P}). In \cite{BV}, Bengochea and Verde, constructed a linear algebraic version of  Mikusi\'nski's theory. With their approach, Bengochea and Verde were able to solve many different equations, such as, certain ode's, difference equations, and fractional differential equations. In this paper we use the basis of the operational calculus introduced in \cite{BV} to construct a  method of algebraic transforms and we solve some Bessel-type equations of order higher o equal than four.

We present next a brief description of our construction. The Mikusiński's operational calculus requires of a vector space of functions $\mathfrak{F} $ and a convolution product under which  $\mathfrak{F} $ becomes in a field, normally the convolution product is defined for functions which Laplace transform is well defined (or any other transform). In our approach, the field is the space of formal Laurent series $\sum_{k\in\mathbb{Z}}a_kp_k,\;a_k\in\mathbb{C},$ where $\{p_k:k\in\mathbb Z\}$ form a group under the operation $p_k*p_m=p_{k+m}$. The elements of $\mathfrak{F}$ are infinite sums  with a finite number of terms with negative index. The product in $\mathfrak{F}$ is defined extending the multiplication of the group $\{p_k:k\in\mathbb Z\}$. Then, we associate to any element of $\mathfrak{F},$ and to some Bessel type equations an algebraic transform. As in any other transforms method, we solve the transformed equation by elementary means, and after partial fractions decomposition. Finally, we apply inverse transform to find the solution of the Bessel type equation .   

This paper is divided as follows, in section \ref{initia} the algebraic setting and basic definitions of the operators to be used are established. In section \ref{BF22} we introduce the Bessel functions of  order $\nu$, and we set the $p_{k,\nu}$ corresponding to the Bessel equations of  order $\nu$ and its respective operator $L_{\nu}$ which corresponds to the  modified left shift operator. At the end of the section we study the Mikusiński's type operational calculus for the operator $L_{0}=1/tDtD.$ We solve two example of differential equations for illustrate our approach. Finally, in section \ref{mocbnu}, we construct an operational calculus for $L_\nu=1/tDtD-\nu^2/t^2.$ As an example, we solve  through this method  the Plum's fourth order differential equation.

\section{Algebraic setting and preliminary results}
\label{initia}
The brief construction to be described here may be found in detail in \cite[section 2]{BV}. Let $\{p_k:k\in\mathbb{Z}\}$ be an abelian group under multiplication defined by $p_k\ast p_n=p_{k+n}$ for $k,n\in\mathbb{Z}$. 
Let $\mathfrak{F}$ be the set of all the formal series of the form $a=\sum_{k\in\mathbb{Z}}a_kp_k,$ where $a_k\in\mathbb{C},\,k\in\mathbb{Z},$ and, either, all the $a_k$ are equal to zero, or there exist an integer $v(a)$ such that $a_k=0$ whenever $k<v(a)$ and $a_{v(a)}\neq 0$. In the first case we write $a=0$ and define $v(a)=\infty.$ Note that  $\mathfrak{F} $ is a linear space over $\mathbb{C}.$ Moreover $\mathfrak{F} $ is a field with multiplication defined by $ab=c=\sum c_np_n$ for $a,b\in\mathfrak{F},$ where
\begin{eqnarray}
\label{2.1}
c_n=\sum_{v(a)\leqslant k\leqslant n-v(b)}a_kb_{n-k}.
\end{eqnarray}
    To each nonzero series $b$ there corresponds a multiplication map  $a\mapsto ab$ this map is clearly linear and invertible. The multiplication map that corresponds to $p_1$ is called the right shift. We denote this correspondence by
   \begin{eqnarray}
   \label{corrs}
   p_1\approx S.
   \end{eqnarray}
The inverse of $S$ denoted $S^{-1}$ is called the left shift and it is denoted by $p_{-1}\approx S^{-1}$. Other operators to be used in this paper are:
\begin{enumerate}
\item $P_n: \mathfrak{F}\to \langle p_n \rangle$ where $\langle p_n \rangle$ is the space generated by $p_n$ and defined by $P a=a_np_n$.
\item $L$ the so called modified left shift which is defined by $La=S^{-1}(I-P_0)a,$ where $I$ is the identity operator. Observe that $Lp_k=p_{k-1}$ for $k\neq 0$ and $Lp_0=0$ and $L$ is not invertible since its kernel is $\langle p_0\rangle.$
\item $A_o:\mathfrak{F}\to\mathbb{C}$ defined by $A_o(a)=a_0.$
\end{enumerate}

\section{Bessel Functions}
\label{BF22}
\subsection{Algebraic Operational Calculus}
\label{AOP}

The Bessel equation of order $\nu$ is given by
\begin{equation}
\label{B1}
R''(t)+\frac{1}{t}R'(t)+\left(1-\frac{\nu^ 2}{t^2}\right) R(t)=0
\end{equation}
which is called the Bessel equation of the first kind of  order $\nu$. To find the solution of (\ref{B1}) using the algebraic operational calculus introduced in \cite{BV}, we define
\begin{eqnarray}
\label{B2}
L_\nu=\frac{1}{t}DtD-\frac{\nu^ 2}{t^ 2}I,\\
\label{Bpk}
p_{k,\nu}=\frac{\left(\frac{t}{2}\right)^{2k+\nu}}{\Gamma(\nu+1+k)k!},
\end{eqnarray}
where for $k\geq 0$ 
$$(-k)!=\frac{(-1)^{k-1}}{(k-1)!}.$$
A simple computation shows 
\begin{equation}
\label{crb}
L_\nu p_{k,\nu}=p_{k-1,\nu},\quad{\rm for}\quad k\neq 0,
\end{equation}
and 
\begin{equation}
\label{crb1}
L_\nu p_{0,\nu}=0,\quad{\rm for}\quad k=0.
\end{equation}
Formulae (\ref{crb}) and (\ref{crb1}) hold for $\nu\in\mathbb{R}$. Observe that the Gamma function 
  $$\Gamma(\nu+1+k)=(\nu+k)!,$$ 
when $\nu$ is a entire.
In this way, $L_\nu$ is the modified left shift introduced in section \ref{initia}.
For instance, the Bessel equation of order zero (\ref{B1}) may be written as follows
\begin{equation}
\label{ocb}
(L_0+I)R=0.
\end{equation}
So after the theory in \cite{BV}, the solution of (\ref{ocb}) (and therefore the solution of (\ref{B1}) with $\nu=0$) is in the space generated by $e_{-1,0;0}=\sum_{k\geq 0}(-1)^k p_{k,0},$ (where $e_{x,0;\nu}$ corresponds to the geometric series $e_{x,0}$ defined after (2.7) p. 334 in \cite{BV}, corresponding to the $L_\nu$) i.e. $\langle e_{-1,0;0}\rangle.$ Since $e_{-1,0;0}$ can be written in terms of $p_{k,0}$ we have that
$$e_{-1,0;0}=\sum_{k\geq 0}(-1)^kp_{k,0}.$$
 We obtain that the solution of the Bessel equation of  order zero in
 terms of this concrete realization  is
\begin{equation}
\nonumber
R(t)=\sum_{k\geq 0}(-1)^k\frac{\left(\frac{t}{2}\right)^{2k}}{(k!)^2}=J_0(t).
\end{equation}
We have obtained by an algebraic purely techniques the well known formula for the Bessel functions of  order zero of the first kind. In general, the Bessel equation of order $\nu$ may be written as
\begin{equation}
\label{ocb1}
(L_\nu+I)R=0.
\end{equation}
We solve (\ref{ocb1})
in the same way that we have solved (\ref{ocb}), so we get the solutions for (\ref{ocb1}) in terms of our series $e_{-1,0;\nu}=J_\nu(t),$ the well known Bessel functions of  order $\nu$.

This method is useful in solving homogeneous differential equations of the form
$$(L_0-c_0I)^{r_0+1}(L_0-c_1)^{r_1+1}\cdots(L_0-c_m)^{r_m+1}y=0,$$
 and even non homogeneous equations of the form
 $$(L_0-c_0I)^{r_0+1}(L_0-c_1)^{r_1+1}\cdots(L_0-c_m)^{r_m+1}y=h,$$
 provided that $h$ is in certain subspace of
 $\mathfrak{F}$  \cite[cf. Corollary 3.2]{BV}.


\subsection{Mikusi\'nski's type Operational Calculus  for $\mathbf{1/tDtD}$}
\label{OCME}
The main object of this paper is to construct an operational calculus over the field $\mathfrak{F}$ for Bessel operators. As an example, in this section we study the operational calculus for $1/tDtD.$ At the end of the section we enunciate the formulae for the operational calculus of
$1/tDtD,$ the operational calculus for the Bessel
 equation of  order zero is constructed without Laplace transforms. 

For  $\nu=0$  in (\ref{Bpk}) we have
$p_{k,0}=(t/2)^{2k}/(k!)^2.$
The set $\{p_{k,0}:k\in\mathbb{Z}\}$ form a group under the operation $p_{m,0}\ast p_{n,0}=p_{m+n,0}.$ Consider the set $\mathfrak{F}=\{\sum a_mp_{m,0}\}$ of the formal Laurent series with  $a_m\in\mathbb{C}.$ Actually $\mathfrak{F}$ is a field with multiplication defined in section
\ref{initia}. In order to obtain a Mikusi\'nski's operational calculus we establish a map between elements $p_k$ and shift operators $S^{k}$. In our case
$p_{1,0}\mapsto S$ and $p_{-1,0}\mapsto S^{-1}$ where $S$ is right shift operator and $S^{-1}$ is left shift operator. Let $S$ correspond to the transformed operator $1/B$ and denote the correspondence with $p_{1,0}$ by
 $$p_{1,0}\approx \frac{1}{B}.$$
 We can write the inverse $S^{-1}$ in terms of the modified left shift  $L_0$ introduced in the last section as
\begin{equation}
\label{MM1}
 S^{-1}=L_0+S^{-1}P_0= B,
\end{equation}
 so that
 $$p_{-1,0}\approx B.$$
 A simple computation shows that $S(L_0+S^{-1}P_0)=(L_0+S^{-1}P_0)S=I.$
 Therefore $$p_{-1,0}\ast p_{1,0}= p_{1,0}\ast p_{-1,0}=p_{0,0}\approx(1/B)B=B/B=I,$$
 where now the symbol $\approx$ means that the equalities in the left side are transformed into the equalities in the right side.
 Provided that $\{p_{k,0}\}$ is a linearly independent set in $\mathfrak{F}$ viewed as linear space over the Complex Numbers, the map $\approx$ is now a linear injective map. We call this map {\it the  Algebraic Transform of the Bessel functions of  order zero} corresponding to $1/tDtD$ and denote it by $\mathcal{A}_{\mathcal{B}_0}$. So that
 $$\mathcal{A}_{\mathcal{B}_0}[p_{k}]=B^{-k},\quad k\in \mathbb{Z}.$$
 In fact, by the discussion in the last section it is easy to see  that if $a\in\mathfrak{F}$ then $p_{-1}\ast a=S^{-1}a\approx 1/tDtD a+4a_0/t^2$ and
 using (\ref{MM1})
 \begin{equation}
\label{imp1}
 \mathcal{A}_{\mathcal{B}_0}\left[\left(\frac{1}{t}DtD+\frac{4A_o}{t^2}\right)a\right]=B\mathcal{A}_{\mathcal{B}_0}[a],
 \end{equation}
  where $A_oa:=a_0.$
 Note that, of course
 $$\mathcal{A}_{\mathcal{B}_0}^{-1}\left[\frac{1}{B^n}\right]=\frac{\left(\frac{t}{2}\right)^{2n}}{(n!)^2},$$
and if $a$ is in the subring $\mathfrak{F}_0=\{a\in\mathfrak{F}:a=\sum_{k\geq 0}a_kp_k\}$
$$\mathcal{A}_{\mathcal{B}_0}^{-1}\left[Ba\right]=\int_0^t\frac{1}{\xi}d\xi\int_0^\xi ua(u)du.$$

\begin{remark}
It is worth to mention that the last expression will never be used in our approach. The main difference of the algebraic operational calculus studied in \cite{BV} and the well known approach of Mikusi\'nski \cite{M} is that in the first one, the integrals and hence the integral transforms and convolutions defined by integrals are not needed at all. In this article the calculations will be purely algebraic, This is the reason for calling $\mathcal{A}_{\mathcal{B}_o}$ algebraic transform.
\end{remark}

\begin{example}
\label{exB1}
{\rm
Derive directly the Bessel function $J_0(\sqrt{\lambda}\,t)$ to obtain}
\begin{eqnarray}
\nonumber
\frac{1}{t}DtDJ_0(\sqrt{\lambda}\,t)&=&\lambda J_0''(\sqrt{\lambda}\,t)+\frac{\lambda}{\sqrt{\lambda}\,t}J_0'(\sqrt{\lambda}\,t)\\\nonumber
&=&-\lambda J_0(\sqrt{\lambda}\,t),
\end{eqnarray}
\end{example}
and respectively
\begin{equation}
\label{3.1p3}
\frac{1}{t}DtDI_0(\sqrt{\lambda}\,t)=\lambda I_0(\sqrt{\lambda}\,t),
\end{equation}
where $I_0$ is the  order zero Bessel function of second kind.
Observe that $A_oI_0(\sqrt{\lambda}\,t)=1,$ since  $I_0(u)=J_0 (iu)=1p_0+i^2p_1+\cdots,$ 
 then
\begin{eqnarray}\label{fI0}
\left(\frac{1}{t}DtD+\frac{4A_o}{t^2}\right)I_0(\sqrt{\lambda}\,t)&=&\lambda I_0(\sqrt{\lambda}\,t)+\frac{4}{t^2}\\\nonumber
&\approx & \\\label{fI02}
B\mathcal{A}_{\mathcal{B}_0}[I_0(\sqrt{\lambda}\,t)]&=&\lambda \mathcal{A}_{\mathcal{B}_0}[I_0(\sqrt{\lambda}\,t)]+BI,
\end{eqnarray}
where $\mathcal{A}_{\mathcal{B}_0}[I_0(\sqrt{\lambda}\,t)]=I+\lambda B^{-1}+\cdots.$
The relation in (\ref{fI02}) follows from \ref{3.1p3} and $4/t^2=p_{-1}\approx B$. So we obtain

\begin{equation}
\label{5.94.1}
(B-\lambda I)\mathcal{A}_{\mathcal{B}_0}[I_0(\sqrt{\lambda}\,t)]=B.
\end{equation}
Therefore
\begin{eqnarray}
\label{5.95}
I_0(\sqrt{\lambda}\,t)=\mathcal{A}_{\mathcal{B}_o}^{-1}\left(\frac{B}{B-\lambda I}\right).
\end{eqnarray}
Further, since $A_oJ_0(\sqrt{\lambda}\,t)=1,$  similar computations as in (\ref{3.1p3}) to (\ref{5.95}) lead to
\begin{eqnarray}
\label{5.97}
J_0(\sqrt{\lambda}\,t)=\mathcal{A}_{\mathcal{B}_o}^{-1}\left(\frac{B}{B+\lambda I}\right).
\end{eqnarray}
Those results are equivalent to the results in \cite[Ch. 5-6]{DR} of Mikusi\'nski's operational calculus for the operator $DtD.$ There, integral convolutions and integrals in general are applied, which are not applied in our approach. Now we can operate formally with the operators in terms of $B$ to obtain the  formulae
\begin{eqnarray}
\label{rusber}
\mbox{ber}(\sqrt{\omega}\, t)&\approx &\frac{B^2}{B^2+\omega^2I},\\
\label{rusbei}
\mbox{bei}(\sqrt{\omega}\, t)&\approx &\frac{\omega B}{B^2+\omega^2I},\\
\frac{1}{2}\left(I_0(\sqrt{\lambda}\,t)+J_0(\sqrt{\lambda}\,t)\right)&\approx &\frac{B^2}{B^2-\lambda ^2I},\\
\frac{1}{2}\left(I_0(\sqrt{\lambda}\,t)-J_0(\sqrt{\lambda}\,t)\right)&\approx &\frac{\lambda B}{B^2-\lambda ^2I},\\\label{TIB0}
\frac{1}{n!}\left(\frac{t^{n}}{2^n\lambda^{n/2}}\right)I_n(\sqrt{\lambda}\,t)&\approx &\frac{B}{(B-\lambda I)^{n+1}},\\
\frac{(-1)^n}{n!}\left(\frac{t^{n}}{2^n\lambda^{n/2}}\right)J_n(\sqrt{\lambda}\,t)&\approx &\frac{B}{(B+\lambda I)^{n+1}}.
\end{eqnarray}
Relations (\ref{rusber}) and (\ref{rusbei}) are obtained by the change of variable $\lambda= i \omega $ in (\ref{5.95}) and the fact that $I_0(\sqrt{i\omega}\,t)=\mbox{ber}\,(\sqrt{\omega}\,t)+\mbox{bei}\,(\sqrt{\omega}\,t).$
Differentiation with respect to $\lambda$ of (\ref{5.95}) and (\ref{5.97}) yields the last two formulae. This formulae are useful when we want to solve ODE's. We finish introducing an expression for initial conditions.

\begin{definition}
\label{gic}
For $a=\sum_{m\geq m_o}a_m p_m\in\mathfrak{F}$ we say that
\begin{eqnarray}
\label{icm}
B^{k}a|_{ic}:=A_o(B^{k}a)= a_k.
\end{eqnarray}
is the generalized $k$-initial condition associated with a problem.
By induction in (\ref{MM1}) we obtain the formula
\begin{equation}
\label{XXIC}
\left(\frac{1}{t}DtD\right)^ku\approx B^ku-a_{k-1}B-a_{k-2}B^2-\cdots-a_0B^k.
\end{equation}
\end{definition}
\hspace{-.6cm}
\begin{remark}
 Note that if $a=\sum_{i\geq 0}a_ip_i,$  then the generalized initial conditions coincide with the initial conditions of an ODE's problem.
\end{remark} 
\begin{example}
Now we are able to solve ordinary differential equations of order higher or equal than two which involve  order zero Bessel functions. As an example we solve the initial conditions problem
\begin{eqnarray}\notag
\label{exmio1}
\left(\frac{1}{t}DtD\right)^2y(t)-(c_1+c_2)\,\frac{1}{t}DtD\;y(t)+c_1c_2\;y(t)=0,\\\notag
By|_{ic}=\alpha,\quad B^2y|_{ic}=\beta.\notag
\end{eqnarray}
After transformations,  using (\ref{XXIC}) and denoting $Y=\mathcal{A}_{\mathcal{B}_0}[y(t)],$ we obtain the equivalent problem
\begin{eqnarray}\notag
\label{exmio1-1}
B^2Y-\beta B-\alpha B^2-(c_1+c_2)BY -(c_1+c_2)\alpha B+c_1c_2Y=0,\\\notag
(B^2-(c_1+c_2)B + c_1c_2I)Y=\alpha B^2+((c_1+c_2)\alpha+\beta) B.
\end{eqnarray}
From the fact that $(B^2-(c_1+c_2)B + c_1c_2I)$ has multiplicative inverse we obtain
\begin{equation}\label{TIB01}
\begin{split}
Y&=\frac{\alpha B^2+((c_1+c_2)\alpha+\beta) B}{(B-c_1I)(B-c_2I)}\\[.3cm]
&=-\frac{(2c_1+c_2)\alpha+\beta}{c_2-c_1}\frac{B}{B-c_1I}+\frac{(c_1+2c_2)\alpha+\beta}{c_2-c_1}\frac{B}{B-c_2I}.
\end{split}
\end{equation}
Finally, we apply the inverse transform in  (\ref{TIB01}) and we achieve
$$y=-\frac{(2c_1+c_2)\alpha+\beta}{c_2-c_1}I_0(\sqrt{c_1}\;t)+\frac{(c_1+2c_2)\alpha+\beta}{c_2-c_1}I_0(\sqrt{c_2}\;t).$$
\end{example}

\section{Mikusi\'nski type Operational Calculus  for $\mathbf{1/tDtD- \nu^2/t^2}$}
\label{mocbnu}
In this section we construct a Mikusi\'nski type operational calculus  over the field $\mathfrak{F}.$  
Recalling the setting in section \ref{AOP} we have
$$L_\nu=\frac{1}{t}DtD-\frac{\nu^2}{t^2}I,$$
and
$$p_{k,\nu}=\frac{\left(\frac{t}{2}\right)^{2k+\nu}}{\Gamma(\nu+1+k)k!}.$$
Now, we proceed as in section \ref{OCME} to construct an operational calculus for the operator $L_\nu$ which is the modified left shift. We make the correspondence
\begin{eqnarray}
p_{1,\nu}&\approx &\frac{1}{B_\nu},
\\
p_{-1,\nu}&\approx & B_\nu,
\end{eqnarray}
  Therefore, we can define a new algebraic transform, which we denote by $\mathcal{A_{B_\nu}}$
and defined by
$$\mathcal{A_{B_\nu}}[p_{k,\nu}]=B^{-k}_{\nu},\quad k\in\mathbb{Z}.$$
Observe that in this concrete realization $S^{-1}P_0a=a_0p_{-1},$ coincides with\break $(t/2)^{-2+\nu}/\Gamma(\nu)A_o a =(t/2)^{-2+\nu}/\Gamma(\nu)a_0$. From the above and the relation (\ref{MM1}) we have for $a\in\mathfrak{F}_\nu$
\begin{eqnarray}
\label{NBNU}
\mathcal{A_{B_\nu}}\left[\left(\frac{1}{t}DtD-\frac{\nu^ 2}{t^2}+\frac{1}{\Gamma(\nu)}\left(\frac{t}{2}\right)^{-2+\nu}A_o \right)a     \right]=B_\nu \mathcal{A_{B_\nu}}[a].
\end{eqnarray}
On the other hand, differentiating directly we obtain
\begin{eqnarray}
\label{bonuj}
\left(\frac{1}{t}DtD-\frac{\nu^2}{t^2}\right)J_\nu(\sqrt{\lambda}\, t)=-\lambda J_\nu(\sqrt{\lambda}\,t),
\end{eqnarray}
and similarly for the modified Bessel functions
\begin{eqnarray}
\label{bonui}
\left(\frac{1}{t}DtD-\frac{\nu^2}{t^2}\right)I_\nu(\sqrt{\lambda}\, t)=\lambda I_\nu(\sqrt{\lambda}\,t).
\end{eqnarray}
Finally, it is easy to verify
\begin{equation}
\label{aoij}
A_o(J_\nu(\sqrt{\lambda}\,t))=A_o(I_\nu(\sqrt{\lambda}\,t))=\lambda^{\nu/2}.
\end{equation}
Taking into account (\ref{NBNU}), (\ref{bonuj}), (\ref{bonui}) and (\ref{aoij}) we obtain the general formulae for the inverse transforms of the Bessel functions $J_\nu,\;I_\nu,$ as in example \ref{exB1}. For instance
\begin{multline}
\label{bonuiff}
\left(\frac{1}{t}DtD-\frac{\nu^2}{t^2}+\frac{1}{\Gamma (\nu)}\left(\frac{t}{2}\right)^{-2+\nu}A_o\right)I_\nu(\sqrt{\lambda}\, t)=\lambda I_\nu(\sqrt{\lambda}\,t)\\+\frac{1}{\Gamma (\nu)}\left(\frac{t}{2}\right)^{-2+\nu}\lambda^{\nu/2},
\end{multline}
and after transformations
\begin{eqnarray}
\label{nt}
B\mathcal{A}_{\mathcal{B}_\nu}[I_\nu(\sqrt{\lambda}\,t)]=\lambda \mathcal{A}_{\mathcal{B}_\nu}[I_\nu(\sqrt{\lambda}\,t)]+\lambda^{\nu/2}B.
\end{eqnarray}
Therefore
\begin{eqnarray}
\label{invij}
I_\nu(\sqrt{\lambda}\,t)&=&\mathcal{A}^{-1}_{\mathcal{B}_\nu}\left(\frac{\lambda^{\nu/2}B_\nu}{B_\nu-\lambda I}\right),\\
J_\nu(\sqrt{\lambda}\,t)&=&\mathcal{A}^{-1}_{\mathcal{B}_\nu}\left(\frac{\lambda^{\nu/2}B_\nu}{B_\nu+\lambda I}\right).
\end{eqnarray}
So we obtain the corresponding Table \ref{T1} of transformations.
\begin{table}
\begin{center}
\begin{tabular}{|l|c|}
\hline
Function $f(t)$ & Transform  $\mathcal{A_{B_\nu}}[f(t)]$\\ \hline
$\displaystyle{\mbox{Ber}_\nu(\sqrt{\omega}\, t)} $ &  $\displaystyle{\frac{B_{\nu}^2}{B_{\nu}^2+\omega^2I}}$\\ \hline
$\displaystyle{\mbox{Bei}_\nu(\sqrt{\omega}\, t)}$ & $\displaystyle{\frac{\omega B_{\nu}}{B_{\nu}^2+\omega^2I}}$\\ \hline
$\displaystyle{\frac{1}{2}\left(I_\nu(\sqrt{\lambda}\,t)+J_\nu(\sqrt{\lambda}\,t)\right)}$ & $\displaystyle{\frac{\lambda^{\nu/2}B_{\nu}^2}{B_{\nu}^2-\lambda ^2I}}$\\ \hline
$\displaystyle{\frac{1}{2}\left(I_\nu(\sqrt{\lambda}\,t)-J_\nu(\sqrt{\lambda}\,t)\right)}$ & $\displaystyle{\frac{\lambda^{\nu/2+1} B_{\nu}}{B_{\nu}^2-\lambda ^2I}}$ \\ \hline 
$\displaystyle{I_{\nu }(\sqrt{\lambda}\,t)}$ & $\displaystyle{\frac{\lambda^{\nu/2}B_{\nu}}{B_{\nu}-\lambda I}}$ \\ \hline
$\displaystyle{J_{\nu }(\sqrt{\lambda}\,t)}$ &  $\displaystyle{\frac{\lambda^{\nu/2}B_{\nu}}{B_{\nu}+\lambda I}}$\\ \hline
\end{tabular}
\label{T1}.
\end{center}
\begin{center}{\small Table 1. Algebraic transforms for $B_{\nu}$}\end{center}
\end{table}
To solve ODE's, we also require generalized initial conditions for $B_\nu$
given by
\begin{eqnarray}
\label{icmnu}
A_o(B_\nu^{k+1}a)&:=&B_\nu^{k+1}a|_{ic}= a_k,
\end{eqnarray}
and the relation
\begin{equation}
\label{cibnu}
\left(\frac{1}{t}DtD-\frac{\nu^2}{t^2}\right)^ku\approx B^k_\nu u-a_{k-1}B_\nu-a_{k-2}B_\nu^2-\cdots-a_0B^k_\nu.
\end{equation}
As an application we solve an equation of fourth order with initial conditions $A_o(B_\nu)=\alpha,\,A_o(B^2_\nu)=\beta$ for $\nu=2.$\\\\
\begin{example}[Solution of Bessel type equations of order higher or equal to 4]
A non trivial application of our operational calculus  is the implementation of a method to find solutions of the  fourth order equation
\begin{equation}
\label{e31}
(ty'')''-((9t^{-1}+8M^{-1}t)y')'=\Lambda t y,
\end{equation}
which appears after separating variables in the solution of the so called Plum equation
$$\Delta^2 u-\gamma\Delta u-\frac{4\gamma}{r^2}u=\Lambda u,$$
where $\Delta$ is the laplacian operator in polar coordinates, \cite[sec. 18]{E}. 
Equation (\ref{e31}) and its relation with the so called Bessel type functions had been profusely studied by Everitt et al. in \cite{EM}, \cite{ESH}, \cite{E}.
Notice that dividing by $t$, and introducing the spectral parameter $\Lambda=\lambda^2(\lambda^2+8/M),$ which is relevant in the study of  Everitt et al., the equation \ref{e31} can be written as
\begin{equation}
\label{e31''}
D^4+\frac{2}{t}D^3-\left(\frac{9}{t^2}+\frac{8}{M}\right)D^2+\left(\frac{9}{t^3}-\frac{8}{Mt}\right)D-\lambda^2\left(\lambda^2+\frac{8}{M}\right) =0.
\end{equation}
Let $L$ be the operator in the left hand side of equation (\ref{e31''}). Then $L$ is the Least Common Left Multiple $(LCLM)$, as defined in \cite{VS}[p. 38]\footnote{Recent versions of some private Mathematics software  had implemented the command DFactorLCLM(L) in the corresponding ODE package to obtain the LCLM factorization of a given operator $L$.}, of the operators
\begin{multline}
\label{L1n}
 L_1=D^2+\frac{\lambda^4 M^2 t^2+8 t^2 \lambda^2 M+16 t^2-48 M }{t (\lambda^4 M^2 t^2+8 t^2 \lambda^2 M+16 t^2-16 M)}D\\
+\frac{\lambda^2 (-4 \lambda^2 M^2-32 M+8 t^2 \lambda^2 M+16 t^2+\lambda^4 M^2 t^2)}{\lambda^4 M^2 t^2+8 t^2 \lambda^2 M+16 t^2-16 M},
\end{multline}
and
\begin{multline}
\label{L2n}
L_2=D^2+\frac{\lambda^4 M^2 t^2+8 t^2 \lambda^2 M+16 t^2-48 M}{t (\lambda^4 M^2 t^2+8 t^2 \lambda^2 M+16 t^2-16 M)}D\\
-\frac{4 M^3 \lambda^4+32 \lambda^2 M^2+16 \lambda^4 M^2 t^2+80 t^2 \lambda^2 M+128 t^2+t^2 \lambda^6 M^3}{M (\lambda^4 M^2 t^2+8 t^2 \lambda^2 M+16 t^2-16 M)}.
\end{multline}
We denote by $L=LCLM(L_1,L_2),$ the fact that $L$ is the $LCLM$ of  $L_1$ and $L_2$. Notice that $L=LCLM(L_1,L_2),$ means that $L$ is rational function times the product  $L_1L_2$, and that the space generated by the solutions of $L_1$ and $L_2$ is the same that the space generated by the solutions of $L.$
If we apply  
the gauge transformation $y(t)\mapsto (4(t^2\lambda^4M+4t^2\lambda^2-16)y(t))/t^2-(32y'(t))/t,$  to the Bessel operator
$$D^2+\frac{1}{t}D-\frac{t^2\lambda^2M+8t^2+4M}{Mt^2},$$
we obtain $L_1.$ Similarly, 
if we apply the gauge transformation $y(t)\mapsto -(4(t^2\lambda^4M^2+12t^2\lambda^2M-16M+32t^2)y(t))/(Mt^2)+(y'(t)32)/t,$ to the Bessel operator
$$D^2+\frac{1}{t}D+\frac{t^2\lambda^2-4}{t^2},$$
we obtain $L_2$ (algorithms to verify the last transformations hold had been implemented in Debeerst Master's thesis \cite{De}).
 Therefore, the space generated by the solutions of (\ref{e31''}) is the same space generated by the solutions of
\begin{equation}
\label{alfin}
\left(D^2+\frac{1}{t}D+\frac{t^2\lambda^2-4}{t^2}\right)\left(D^2+\frac{1}{t}D-\frac{t^2\lambda^2M+8t^2+4M}{Mt^2}\right)y(t)=0,
\end{equation}
or equivalently 
\begin{equation}
\label{alfin1}
\left(\frac{1}{t}DtD-\frac{4}{t^2}+\lambda^2\right)\left(\frac{1}{t}DtD-\frac{4}{t^2}-\left(\lambda^2+\frac{8}{M}\right)\right)y(t)=0.
\end{equation}
After expanding we obtain
\begin{equation}
\label{alfin11}
{\textstyle\left[\left(\frac{1}{t}DtD-\frac{4}{t^2}\right)^2-\left(\left(\lambda^2+\frac{8}{M}\right)-\lambda^2\right)\left(\frac{1}{t}DtD-\frac{4}{t^2}\right)-\left(\lambda^2+\frac{8}{M}\right)\lambda^2\right]y(t)=0.}
\end{equation}
From (\ref{icmnu}) and (\ref{cibnu}), with $\nu=2,$ we obtain,
after factoring
\begin{eqnarray}
\label{BY}
(B_2+\lambda^2I)(B_2- (\lambda^2+8/M)I)Y=\alpha B_2^2+((\lambda^2+8/M)-\lambda^2)\alpha+\beta) B_2,
\end{eqnarray}
where $A_o(B_2)=\alpha,\,A_o(B^2_2)=\beta.$
So as in Example 2, after partial fraction decomposition, and  taking inverse transforms from Table 1, we obtain
$$y=-\frac{(-\lambda^2+8/M)\alpha+\beta}{\lambda(2\lambda^2+8/M)}J_2(\lambda\;t)+\frac{(\lambda^2+16/M)\alpha+\beta}{(\lambda^2+8/M)(2\lambda^2+8/M)}I_2(\sqrt{\lambda^2+8/M }\;t).$$
Finally, from the well known formulae
\begin{eqnarray*}
J_2(\lambda t)&=&-J_0(\lambda t)+\frac{2}{\lambda t}J_1(\lambda t),\\
I_2(t\sqrt{\lambda^2+8/M} )&=& I_0(t\sqrt{\lambda^2+8/M} )-\frac{2}{t\sqrt{\lambda^2+8/M}} I_1(t\sqrt{\lambda^2+8/M} ).
\end{eqnarray*}
we obtain the same space $S$ as in \cite{E}, generated by solutions of (\ref{e31''})  uniformly bounded in $[0,\infty),$   i.e.
$$S=\left\langle J_0(\lambda t),\frac{J_1(\lambda t)}{\lambda t},I_0(t\sqrt{\lambda^2+8/M}),\frac{I_1(t\sqrt{\lambda^2+8/M})}{t\sqrt{\lambda^2+8/M}} \right\rangle,$$
where $\lambda\in \mathbb{C},M\in (0,\infty)$, $\Lambda=-\lambda^2(\lambda^2+8/M).$
\end{example}

\section{Conclusions}
We had develop a transforms method to find bounded solutions of Bessel type equations of arbitrary  order $\nu$ in terms of first kind Bessel functions and modified first kind Bessel functions. Nevertheless our examples consist only in homogeneous equations of order higher or equal than two, our method can be applied to non-homogeneous problems if the non homogeneous term is in the space   $\mathfrak{F},$ defined in section \ref{initia}. 

\subsection{Acknowledgements}
We want to thank to Luis Verde-Star for his careful reading and some corrections.

\end{document}